\documentclass[11pt,reqno]{amsart}
\usepackage{amsfonts,amsmath,amssymb}
\newtheorem{thm}{Theorem}[section]

\newtheorem{lemma}[thm]{Lemma}

\usepackage{graphicx}
\usepackage{epstopdf}

\title[Jacobi--Trudy formula for generalised Schur polynomials] {Jacobi--Trudy formula for generalised Schur polynomials}
\author{A.N. Sergeev}\address{Department of Mathematical Sciences,
Loughborough University, Loughborough LE11 3TU, UK}
\email{A.N.Sergeev@lboro.ac.uk}

\author{A.P. Veselov}
\address{Department of Mathematical Sciences,
Loughborough University, Loughborough LE11 3TU, UK and Department of Mathematics and Mechanics, Moscow State University, Moscow, 119899, Russia}
\email{A.P.Veselov@lboro.ac.uk}

\begin{document}
\begin{abstract} Jacobi--Trudy formula for a generalisation of Schur polynomials related to any sequence of orthogonal polynomials in one variable is given. As a corollary we have Giambelli formula for generalised Schur polynomials.
\end{abstract}

\maketitle

\section{Introduction}
The classical Jacobi--Trudy formula expresses the Schur polynomials 
\begin{equation}
\label{JT0clas}
S_{\lambda} (x_{1},\dots,x_{n})=\frac{\left|\begin{array}{cccc}
   x_{1}^{\lambda_{1}+n-1}& x_{2}^{\lambda_{1}+n-1}& \ldots&x_n^{\lambda_{1}+n-1}\\
  x_{1}^{\lambda_{2}+n-2}&x_2^{\lambda_{2}+n-2}& \ldots &x_n^{\lambda_{2}+n-2}\\
\vdots&\vdots&\ddots&\vdots\\
 x_1^{\lambda_{n}}&x_2^{\lambda_{n}}& \ldots&x_n^{\lambda_{n}}\\
 \end{array}\right|}{\Delta_{n}(x)}
 \end{equation}
 where $\lambda=(\lambda_1,\dots, \lambda_n)$ is a partition and $\Delta_{n}(x)=\prod_{i<j}^n(x_{i}-x_{j}),$ as the determinant
\begin{equation}
\label{JTclas}
S_{\lambda}(x_{1},\dots,x_{n})=\left|\begin{array}{cccc}
   h_{\lambda_{1}}&h_{\lambda_{1}+1}& \ldots &h_{\lambda_{1}+l-1}\\
  h_{\lambda_{2}-1}&h_{\lambda_{2}}& \ldots &h_{\lambda_{2}+l-2}\\
\vdots&\vdots&\ddots&\vdots\\
  h_{\lambda_{l}-l+1}&h_{\lambda_{l}-l+2}& \ldots &h_{\lambda_{l}}\\
 \end{array}\right|
 \end{equation}
where $l=l(\lambda)$ and $h_i=h_i(x_1,\dots,x_n)$ are the complete symmetric polynomials (see \cite{Mac}). Note that these polynomials $h_i$ are particular case of Schur polynomials $S_{\lambda},$ corresponding to the partition $\lambda=(i)$ consisting of one part.

In this note we give a version of this formula, which is valid for the following generalisation of Schur polynomials related to any sequence of orthogonal  polynomials in one variable. 

More precisely, 
let  $\{\varphi_{i}(z)\},\; i=0,1,2,\dots $ be a sequence of polynomials in one variable,
which satisfy a three-term recurrence relation 
\begin{equation}
\label{relaph}
z\varphi_{i}(z)=\varphi_{i+1}(z)+a(i)\varphi_{i}(z)+b(i)\varphi_{i-1}(z)
\end{equation}
with $\varphi_0\equiv 1, \,\varphi_{-1}\equiv 0$ (for example, a sequence of the orthogonal polynomials \cite{Sz}). 
The corresponding {\it generalised Schur polynomials} $S(x_1,\dots,x_n|a,b)$ are defined for any partition $\lambda$ and two infinite sequences $a=\{a_i\},\, b=\{b_i\}$ by the Weyl-type formula
\begin{equation}
\label{JT0}
S_{\lambda}(x_{1},\dots,x_{n}|a,b)=\frac{\left|\begin{array}{cccc}
   \varphi_{\lambda_{1}+n-1}(x_{1})&\varphi_{\lambda_{1}+n-1}(x_{2})& \ldots &\varphi_{\lambda_{1}+n-1}(x_{n})\\
 \varphi_{\lambda_{2}+n-2}(x_{1})&\varphi_{\lambda_{2}+n-2}(x_{2})& \ldots &\varphi_{\lambda_{2}+n-2}(x_{n})\\
\vdots&\vdots&\ddots&\vdots\\
  \varphi_{\lambda_{n}}(x_{1})&\varphi_{\lambda_{n}}(x_{2})& \ldots &\varphi_{\lambda_{n}}(x_{n})\\
 \end{array}\right|}{\Delta_{n}(u)}.
 \end{equation}
For $n=1$ we assume that $\Delta_1\equiv 1,$ so for $\lambda=(i)$ we have $S_{\lambda}(x_{1}|a,b)= \varphi_i(x_1).$ If the initial sequence $\varphi_i(z)$ was orthogonal with measure $d\mu(z)$ the polynomials $S_{\lambda}(x_{1},\dots,x_{n}|a,b)$ are orthogonal with respect to the measure
\begin{equation}
\label{mes}
\Omega(z)= \Delta^2_{n}(z) \prod_{i=1}^n d\mu(z_i).
 \end{equation}

Alternatively, the generalised Schur polynomials can be defined in this case as the polynomials of the triangular form
\begin{equation}
\label{trian}
S_{\lambda}(x_{1},\dots,x_{n}|a,b)=\sum_{\mu\preceq \lambda} K_{\lambda,\mu}(a,b)m_{\mu}(x_{1},\dots,x_{n}),
\end{equation}
which are orthogonal with respect to the measure (\ref{mes}).
Here $m_{\mu}(x_{1},\dots,x_{n})$ are monomial symmetric polynomials \cite{Mac} and the notation $\mu\preceq \lambda$ means that $\mu_1+\dots +\mu_k \leq \lambda_1+\dots +\lambda_k$ for all $1 \leq k\leq n.$ 

When $\varphi_i(z)$ is the sequence of classical Jacobi polynomials \cite{Sz}, the corresponding generalised Schur polynomials coincide with the multidimensional Jacobi polynomials with parameter $\theta=1$ (see Lassalle \cite{Lass} and  Okounkov-Olshanski \cite{OO}).

Denote the polynomials $S_{\lambda}(x_{1},\dots,x_{n}|a,b)$ with $\lambda=(i,0,\dots,0)$ as $h_i(x)$ and extend this sequence for negative $i$ by assuming that $h_i(x)\equiv 0$ for $ i<0.$
Extend also the sequence of coefficients $a(i)$ and $b(i)$ to the negative $i$ arbitrarily and define recursively the polynomials $h_i^{(r)}(x_1,\dots,x_n)$ by
the relation
\begin{equation}
\label{relah}
h_i^{(r+1)}=h_{i+1}^{(r)}+ a(i+n-1)h_i^{(r)}+b(i+n-1)h_{i-1}^{(r)}
\end{equation}
with initial data $ h_i^{(0)}(x)= h_i(x).$ One can check that 
$h_i^{(r)}(x)\equiv 0$ whenever $i+r<0$ and
that the definition of the polynomials $h_i^{(r)}(x_1,\dots,x_n)$ does not depend on the extension of the coefficients to the negative $i$ provided
\begin{equation}
\label{relar}
r\leq i+2n-2.
\end{equation}
 In particular, all the entries of the formula (\ref{JT}) below are well defined.

Our main result is the following

\begin{thm} 
\label{JTthm}
The generalised Schur polynomials satisfy the following Jacobi-Trudy formula:
\begin{equation}
\label{JT}
S_{\lambda}(x_{1},\dots,x_{n}|a,b)=\left|\begin{array}{cccc}
   h_{\lambda_{1}}&h^{(1)}_{\lambda_{1}}& \ldots &h^{(l-1)}_{\lambda_{1}}\\
  h_{\lambda_{2}-1}&h^{(1)}_{\lambda_{2}-1}& \ldots &h^{(l-1)}_{\lambda_{2}-1}\\
\vdots&\vdots&\ddots&\vdots\\
  h_{\lambda_{l}-l+1}&h^{(1)}_{\lambda_{l}-l+1}& \ldots &h^{(l-1)}_{\lambda_{l}-l+1}\\
 \end{array}\right|
 \end{equation}
where $l=l(\lambda)$. 
\end{thm}

This gives a universal proof of the Jacobi--Trudy and Giambelli formulas  for usual Schur polynomals as well as for the characters of symplectic and orthogonal Lie algebras (see \cite{FH}) and for the factorial Schur polynomials \cite{Mac,BL}. Another interesting case, which seems to be new, is the Jacobi--Trudy formula for the multidimensional Jacobi polynomials with parameter $\theta=1.$ 

\section{Proof}

We start with the following lemma.
\begin{lemma} The following equality holds
\begin{equation}
\label{fl}
h_i^{(r)}(x_1,\dots,x_n)-x_1h_i^{(r-1)}(x_1,\dots,x_n)=h_{i+1}^{(r-1)}(x_2,\dots,x_n)
\end{equation}
for all $r,i$ satisfying the relation (\ref{relar}).
\end{lemma}

\begin{proof} The proof is by induction in $r.$ When $r=1$ we have from definition 
$$
h_i^{(1)}(x_1,\dots,x_n)-x_1h_i^{(0)}(x_1,\dots,x_n)= h_{i+1}(x_{1},\dots,x_{n})+a(i+n-1)h_{i}(x_{1},\dots,x_{n})
$$
$$
+b(i+n-1)h_{i-1}(x_{1},\dots,x_{n})-x_{1}h_{i}(x_{1},\dots,x_{n})=
$$
$$
=\Delta_{n}(x)^{-1}\left|\begin{array}{cccc}
   0&(x_{2}-x_{1})\varphi_{i+n-1}(x_{2})& \ldots &(x_{n}-x_{1})\varphi_{i+n-1}(x_{n})\\
 \varphi_{n-2}(x_{1})&\varphi_{n-2}(x_{2})& \ldots &\varphi_{n-2}(x_{n})\\
\vdots&\vdots&\ddots&\vdots\\
 1&1& \ldots &1\\
 \end{array}\right|.
 $$
Subtracting the first column from the others we get
$$
=\Delta_{n}(x)^{-1}\left|\begin{array}{cccc}
   0&(x_{2}-x_{1})\varphi_{i+n-1}(x_{2})& \ldots &(x_{n}-x_{1})\varphi_{l+n-1}(x_{n})\\
 \varphi_{n-2}(x_{1})&\varphi_{n-2}(x_{2})-\varphi_{n-2}(x_{1})& \ldots &\varphi_{n-2}(x_{n})-\varphi_{n-2}(x_{1})\\
\vdots&\vdots&\ddots&\vdots\\
 1&0& \ldots &0\\
 \end{array}\right|=
 $$
$$
=\Delta_{n-1}(x)^{-1}\left|\begin{array}{ccc}
  \varphi_{i+n-1}(x_{2})& \ldots &\varphi_{i+n-1}(x_{n})\\
 \frac{\varphi_{n-2}(x_{2})-\varphi_{n-2}(x_{1})}{x_{2}-x_{1}}& \ldots & \frac{\varphi_{n-2}(x_{n})-\varphi_{n-2}(x_{1})}{x_{n}-x_{1}}\\
\vdots&\ddots&\vdots\\
  \frac{\varphi_{1}(x_{2})-\varphi_{1}(x_{1})}{x_{2}-x_{1}}&\ldots&  \frac{\varphi_{1}(x_{n})-\varphi_{1}(x_{1})}{x_{n}-x_{1}}\\
 \end{array}\right|=h_{i+1}(x_{2},\dots,x_{n}).
 $$
 The induction step is straightforward check using the relation (\ref{relah}).
 \end{proof}

Now we are ready to prove the Jacobi--Trudy  formula. The proof is by  induction in $l=l(\lambda)$. 
If $l=1$ then the formula follows from the definition of $h_i.$ Suppose that $l>1$. We will use the bracket $\{g(x_1,\dots,x_n)\}$ to denote the result of the alternation:
$$\{g(x_1,\dots,x_n)\}=\sum_{w\in S_{n}}\varepsilon(w)g(x_{w(1)},\dots,x_{w(n)}).$$
We claim that
$$
\{h_{i}^{(r)}(x_1,\dots,x_n) x_{1}^{n-1}x_{2}^{n-2}\dots x_{n}^0\}=\{x_1^r \varphi_{i+n-1}(x_{1})x_{2}^{n-2}\dots x_{n}^0\}
$$
for any $r\le i+2n-2.$ Indeed, for $r=0$ this true by definition and the induction step follows easily from relations (\ref{relaph}) and (\ref{relah}).
From this we have
$$
\left\{\left|\begin{array}{cccc}
   h_{\lambda_{1}}&h^{(1)}_{\lambda_{1}}& \ldots &h^{(l-1)}_{\lambda_{1}}\\
  h_{\lambda_{2}-1}&h^{(1)}_{\lambda_{2}-1}& \ldots &h^{(l-1)}_{\lambda_{2}-1}\\
\vdots&\vdots&\ddots&\vdots\\
  h_{\lambda_{l}-l+1}&h^{(1)}_{\lambda_{l}-l+1}& \ldots &h^{(l-1)}_{\lambda_{l}-l+1}\\
 \end{array}\right|
x_{1}^{n-1}x_{2}^{n-2}\dots x_{n}^0\right\}
 $$
$$
=\left\{\left|\begin{array}{cccc}
   \varphi_{\lambda_{1}+n-1}(x_1)&x_1 \varphi_{\lambda_{1}+n-1}(x_1)& \ldots &x_1^{l-1}\varphi_{\lambda_{1}+n-1}(x_1)\\
  h_{\lambda_{2}-1}&h^{(1)}_{\lambda_{2}-1}& \ldots &h^{(l-1)}_{\lambda_{2}-1}\\
\vdots&\vdots&\ddots&\vdots\\
  h_{\lambda_{l}-l+1}&h^{(1)}_{\lambda_{l}-l+1}& \ldots &h^{(l-1)}_{\lambda_{l}-l+1}\\
 \end{array}\right|
x_{2}^{n-2}\dots x_{n}^0\right\}.
 $$
Multiply every column except  the last one by $x_{1}$ and subtract it from the next column. Then by lemma this expression is equal to
$$
\left\{\left|\begin{array}{cccc}
   \varphi_{\lambda_{1}+n-1}(x_1)&0& \ldots &0\\
  h_{\lambda_{2}-1}&\hat h_{\lambda_{2}-1}& \ldots &\hat h^{(l-2)}_{\lambda_{2}-1}\\
\vdots&\vdots&\ddots&\vdots\\
  h_{\lambda_{l}-l+1}&\hat h_{\lambda_{l}-l+1}& \ldots &\hat h^{(l-2)}_{\lambda_{l}-l+1}\\
 \end{array}\right|
x_{2}^{n-2}\dots x_{n}^0\right\}$$
$$
=\left\{\varphi_{\lambda_{1}+n-1}(x_1)\left |\begin{array}{ccc}
   \hat h_{\lambda_{2}-1}& \ldots &\hat h^{(l-2)}_{\lambda_{2}-1}\\
\vdots&\ddots&\vdots\\
 \hat h_{\lambda_{l}-l+1}& \ldots &\hat h^{(l-2)}_{\lambda_{l}-l+1}\\
 \end{array}\right|
x_{2}^{n-2}\dots x_{n}^0\right\},$$
where $\hat h^{(r)}_{i}=h^{(r)}_{i}(x_{2},\dots, x_{n}).$ By induction this is equal to
$$\{\varphi_{\lambda_{1}+n-1}(x_1)\varphi_{\lambda_{2}+n-2}(x_2)\dots\varphi_{\lambda_{n}}(x_n)\},
$$
which by definition coincides with $f_{\lambda}.$ 
This completes the proof of the main theorem.

\section{Giambelli formula}

As a corollary we have the  following Giambelli  formula for generalised Schur functions. \footnote{ We are very grateful to G. Olshanski, who pointed out this to us.}
Let us denote the generalised Schur polynomials  corresponding to the hook Young diagrams as $$S_{(u|v)}(x)=S_{(u+1,1^{(v)})}(x).$$

\begin{thm} 
The generalised Schur polynomials satisfy the following Giambelli formula
\begin{equation}
\label{Gia}
S_{\lambda}(x_{1},\dots,x_{n}|a,b)=\left|\begin{array}{cccc}
   S_{(\lambda_{1}-1|\lambda_{1}^{'}-1)}&S_{(\lambda_{1}-1|\lambda_{2}^{'}-2)}& \ldots &S_{(\lambda_{1}-1|\lambda_{r}^{'}-r)}\\
  S_{(\lambda_{2}-2|\lambda_{1}^{'}-1)}&S_{(\lambda_{2}-2|\lambda_{2}^{'}-2)}& \ldots &S_{(\lambda_{2}-2|\lambda_{r}^{'}-r)} \\
\vdots&\vdots&\ddots&\vdots\\
   S_{(\lambda_{r}-r|\lambda_{1}^{'}-1)}&S_{(\lambda_{r}-r|\lambda_{2}^{'}-2)}& \ldots &S_{(\lambda_{r}-r|\lambda_{r}^{'}-r)}\\
 \end{array}\right|
 \end{equation}
where $r$ is the number of the diagonal boxes of $\lambda$.
\end{thm}
\begin{proof} The proof follows the same line as Macdonald's proof of the usual Giambelli formula (see \cite{Mac}, Ch.1, Section 3, Example 21), but we give the proof here for  the reader's convenience. 

From the Theorem \ref{JTthm} we see that 
$$
S_{(u|v)}(x)=\left|\begin{array}{cccc}
   h_{u+1}&h^{(1)}_{u+1}& \ldots &h^{(v)}_{u+1}\\
  h_{0}&h^{(1)}_{0}& \ldots &h^{(v)}_{0}\\
\vdots&\vdots&\ddots&\vdots\\
  h_{1-v}&h^{(1)}_{1-v}& \ldots &h^{(v)}_{1-v}\\
 \end{array}\right|
 $$ 
 In this formula $u \ge 0,\: v\ge0$, but we can define the functions $S_{(u|v)}(x)$ by the same formula for all  integers $u$ and nonnegative integers $v.$ It is easy to check that  this defines them correctly and that for $u$ negative $S_{(u|v)}(x)=0$ except when $u+v=-1$, in which case $S_{(u|v)}(x)=(-1)^v.$
 
 Now consider the following matrix of the size $j\times (j+1)$ 
$$
 H^{(j)}=\left(\begin{array}{cccc}
   h_{0}&h_{0}^{(1)}& \ldots &h_{0}^{(j)}\\
  h_{-1}&h_{-1}^{(1)}& \ldots &h_{-1}^{(j)}\\
\vdots&\vdots&\ddots&\vdots\\
 h_{1-j}&h_{1-j}^{(1)}& \ldots &h_{1-j}^{(j)}\\
 \end{array}\right)
$$
and  denote by $\Delta_{i}^{(j)},\, 1\le i\le j+1$ the determinant  of its  sub-matrix without the $i$-th column multiplied by $(-1)^{i-1}.$  If $i>j+1$ we  set  by definition $\Delta_{i}^{(j)}=0$. One can check also that $\Delta_k^{(k-1)}=(-1)^{k-1}.$

For any partition $\lambda$ consider the matrices 
$$
A=\left(\begin{array}{cccc}
   h_{\lambda_1}&h_{\lambda_1}^{(1)}& \ldots &h_{\lambda_1}^{(l-1)}\\
  h_{\lambda_2-1}&h_{\lambda_2-1}^{(1)}& \ldots &h_{\lambda_2-1}^{(l-1)}\\
\vdots&\vdots&\ddots&\vdots\\
 h_{\lambda_l-l+1}&h_{\lambda_l-l+1}^{(1)}& \ldots &h_{\lambda_l-l+1}^{(l-1)}\\
 \end{array}\right), \,\,
B=\left(\begin{array}{cccc}
   \Delta_{1}^{(l-1)}&\Delta_{1}^{(l-2)}& \ldots &\Delta_{1}^{(0)}\\
  \Delta_{2}^{(l-1)}&\Delta_{2}^{(l-2)}& \ldots &\Delta_{2}^{(0)}\\
\vdots&\vdots&\ddots&\vdots\\
\Delta_{l}^{(l-1)}& \Delta_{l}^{(l-2)}& \ldots & \Delta_{l}^{(0)}\\
 \end{array}\right).
$$
Note that $B$ is upper-triangular with respect to the anti-diagonal with the anti-diagonal elements $(-1)^{k-1}$, so the determinant of $B$ is identically equal to 1, while the determinant of $A$ by Theorem \ref{JTthm} coincides with $S_{\lambda}(x_{1},\dots,x_{n}|a,b).$ From linear algebra and definition of $S_{(u|v)}$ we have
$$
AB=\left(\begin{array}{cccc}
   S_{(\lambda_{1}-1|l-1)}&S_{(\lambda_{1}-1|l-2)}& \ldots &S_{(\lambda_{1}-1|0)}\\
  S_{(\lambda_{2}-2|l-1)}&S_{(\lambda_{2}-2|l-2)}& \ldots &S_{(\lambda_{2}-2|0)} \\
\vdots&\vdots&\ddots&\vdots\\
   S_{(\lambda_{l}-l|l-1)}&S_{(\lambda_{l}-l|l-2)}& \ldots &S_{(\lambda_{l}-l|0)}\\
 \end{array}\right).
$$
Taking the determinants of both sides we see that $S_{\lambda}(x_{1},\dots,x_{n}|a,b)$ equals to the determinant of the last matrix. In this matrix there are many zeros since for $k>r$ we have $\lambda_k-k<0$ and therefore $S_{(\lambda_{k}-k|l-j)}=(-1)^{l-j}$ if $\lambda_k-k+l-j=-1$ and $0$ otherwise. This means that in the $k$-th row with $k>r$ there is only one non-zero element $(-1)^{l-j}$ with $l-j=k-\lambda_k-1.$ This reduces the calculation of the determinant to the $r\times r$ matrix with the remaining columns having the numbers $\lambda'_j-j,\,\, j=1,\dots,r.$ Indeed, for any $\lambda$ of length $l$ with $r$ boxes on the diagonal the union of two sets $\{k-\lambda_k-1\}, k=r+1,\dots,l$ and $\{\lambda'_j-j\}, j=1,\dots, r$ is the set
$\{0,1,2, \dots, l-1\}$ as it follows, for example, from the identity
$$
\sum_{i=1}^lt^{i}(1-t^{-\lambda_{i}})=\sum_{j=1}^r(t^{\lambda_{j}^{'}-j+1}-t^{j-\lambda_{j}})
$$
(see \cite{Mac}, Ch.1, Section 1, Example 4). The check of the sign completes the proof.
\end{proof}

\section{Particular cases}

As a corollary we have the following well-known cases of the Jacobi--Trudy formula.

1. When $a(i)=b(i)=0$ for all $i\ge 0$ we have $\varphi_i(z)=z^i$ and (\ref{JT}) clearly coincides with the usual Jacobi--Trudy formula for Schur polynomials.

2. The characters of the orthogonal Lie algebra $so(2n+1)$ correspond to the case when $a(i)=0,\, b(i)=1$ for $i>0$ and $a(0)=-1,\, b(0)=0$ and the polynomials $\varphi_{i}(z)=x^i+x^{i-1}+\dots+x^{-i},\: z=x+x^{-1}.$ Using the recurrence relation (\ref{relah}), having in this case the form
$$h_i^{(r+1)}=h_{i+1}^{(r)}+ h_{i-1}^{(r)}$$ we can rewrite the general Jacobi-Trudy formula (\ref{JT}) 
in the form known in representation theory (see Prop. 24.33 in Fulton-Harris \cite{FH}):
the character $\chi_{\lambda}$ is the determinant of the $l\times l$ matrix whose $i$-th row is
$$(h_{\lambda_i-i+1} \quad h_{\lambda_i-i+2} +h_{\lambda_i-i} \quad h_{\lambda_i-i+3} +h_{\lambda_i-i-1} \,\,  \dots \,\,  h_{\lambda_i-i+l} +h_{\lambda_i+2-l} ).$$

The same is true for the characters of the even orthogonal Lie algebra $so(2n),$
where $a(i)=0$ for all $i\ge 0$ and $b(i)=1$ for $i>1$ with $b(1)=2$ and $\varphi_{i}(z)=x^i+x^{-i},\: z=x+x^{-1}$ (see Prop. 24.44 in \cite{FH}) and for the symplectic Lie algebra $sp(2n),$
when $a(i)=0,\, b(i)=1$ for all $i\ge0$ and $\varphi_{i}(z)=x^i+x^{i-2}+\dots+x^{-i},\: z=x+x^{-1}$ (Prop. 24.22 in \cite{FH}).
Note that the change of $a(0)$ and $b(1)$ does not affect the definition of the relevant $h^{(r)}_i$ for $r>0.$

3. The {\it factorial Schur polynomials} \cite{BL} correspond to the special case when $b_i=0,$ so that
$$\varphi_i(z)=(z-a(0))(z-a(1))\dots (z-a(i-1)),\quad i>0.$$ The Jacobi--Trudy formula for them can be found in \cite{Mac}, Ch.1, Section 3, Example 20.

4. For $a(i),\,b(i)$ given by 
\begin{equation}
\label{bca}
a(x)=-\frac{2p(p+2q+1)}{(2x-p-2q-1)(2x-p-2q+1)},
\end{equation}
\begin{equation}
\label{bcb}
 b(x)=\frac{2x(2x-2q-1)(2x-2p-2q-1)(2x-2p-4q-2)}{(2x-p-2q)(2x-p-2q-1)^2(2x-p-2q-2)}
\end{equation}
 we have the Jacobi-Trudy formula for the multidimensional Jacobi polynomials with $k=-1,$ which seems to be new.

 \section{Infinite-dimensional and super versions}
 
Let us assume now that the coefficients $a(i)$ and $b(i)$ of the recurrence relation are {\it rational functions} of $i.$ In that case we can define the {\it generalised Schur functions} (which are the infinite-dimensional version of $S_{\lambda}(x|a,b)$) in the following way (cf. Okounkov-Olshanski \cite{OO}).
 
 First note that the generalised Schur polynomials (\ref{JT0}) are the linear combination of the usual Schur polynomials
 $$S_{\lambda}(x_1,\dots,x_n|a,b)= \sum_{\mu\subseteq \lambda} c_{\lambda,\mu}(n|a,b)S_{\mu}(x_1,\dots,x_n),$$
where $c_{\lambda,\mu}(n|a,b)$ are some rational functions of $n.$ The generalised Schur functions depend on the additional parameter $d$ and defined by
\begin{equation}
\label{gsch}
S_{\lambda}(x|a,b;d)= \sum_{\mu\subseteq \lambda} c_{\lambda,\mu}(d|a,b)S_{\mu}(x),
\end{equation}
where $S_{\mu}(x)$ are the usual Schur functions \cite{Mac}. They satisfy the following
 infinite dimensional version of the Jacobi--Trudy formula:
\begin{equation}
\label{JTinfi}
S_{\lambda}(x|a,b;d)=\left|\begin{array}{cccc}
   h_{\lambda_{1}}&h^{(1)}_{\lambda_{1}}& \ldots &h^{(l-1)}_{\lambda_{1}}\\
  h_{\lambda_{2}-1}&h^{(1)}_{\lambda_{2}-1}& \ldots &h^{(l-1)}_{\lambda_{2}-1}\\
\vdots&\vdots&\ddots&\vdots\\
  h_{\lambda_{l}-l+1}&h^{(1)}_{\lambda_{l}-l+1}& \ldots &h^{(l-1)}_{\lambda_{l}-l+1}\\
 \end{array}\right|.
 \end{equation}
 Here $l=l(\lambda)$, $h_i=S_{\lambda}(x|a,b;d)$ with $\lambda=(i)$ if $i\geq 0$ and $h_i\equiv 0$ if $i<0,$ $h_i^{(r)}= h_i^{(r)}(x|a,b;d)$ are defined for generic $d$ by the recurrence relation
\begin{equation}
\label{relahinf}
h_i^{(r+1)}=h_{i+1}^{(r)}+ a(i+d-1)h_i^{(r)}+b(i+d-1)h_{i-1}^{(r)}
\end{equation}
with initial data $ h_i^{(0)}= h_i.$

The {\it generalised super Schur polynomials} $S_{\mu}(x_1,\dots,x_n; y_1,\dots,y_m|a,b)$ can be defined by the same formula (\ref{gsch}), where the Schur functions should be replaced by the super Schur polynomials $S_{\mu}(x_1,\dots,x_n; y_1,\dots,y_m)$ (see e.g. \cite{Mac}) and $d$ must be specialised as the superdimension $d=n-m$ (provided the coefficients have no poles at $d=n-m$). Alternatively, $S_{\mu}(x_1,\dots,x_n; y_1,\dots,y_m|a,b)$ is the image of the corresponding generalised Schur function $S_{\mu}(x|a,b;d)$ under the homomorphism $\phi$ sending the power sums $p_k \in \Lambda$ to the super power sums $x_1^k+\dots+x_n^k-y_1^k-\dots-y_m^k$ with $d=n-m.$ 
In the case of factorial Schur polynomials their super version had been introduced in a different way by Molev \cite{Molev}.

An important  example corresponds to the sequences (\ref{bca}),(\ref{bcb}). In this case the generalised Schur functions coincide with the Jacobi symmetric functions with parameter $k=-1$ (see \cite{SV1}).
These functions and their super versions play an important role in representation theory of the orthosymplectic Lie superalgebras
\cite{SV2}.

Finally we would like to mention that the Jacobi-Trudy formula (\ref{JT}) can be rewritten in a dual form in terms of the conjugate partition in the spirit of  Macdonald \cite{Mac} (Ch.1, Section 3, Example 21) and Okounkov-Olshanski \cite{OO2}, Section 13.

\section{Acknowledgements}
We are grateful to A. Okounkov and G. Olshanski for useful comments.

This work has been partially supported by EPSRC (grant EP/E004008/1) and by the
European Union through the FP6 Marie Curie RTN ENIGMA (contract
number MRTN-CT-2004-5652) and through ESF programme MISGAM.

\end{document}